\documentclass[12pt]{amsart}
\usepackage{graphicx}
\usepackage{amssymb}
\usepackage{amsfonts}
\usepackage{amsmath}
\usepackage{array}
\usepackage{rotating}

\newtheorem{example}{Example}[section]
\newtheorem{note}[example]{Note}
\newtheorem{theorem}[example]{Theorem}

\newtheorem{definition}[example]{Definition}
\newtheorem{proposition}[example]{Proposition}

\def\bp{BP}
\def\match{Ma}

\def\FQSym{{\bf FQSym}}
\def\MQSym{{\bf MQSym}}
\def\PQSym{{\bf PQSym}}

\def\NCSF{{\bf Sym}}
\def\QSym{{\it QSym}}
\def\FSym{{\bf FSym}}
\def\PBT{{\bf PBT}}
\def\sym{{\it Sym}}

\def\Std{{\rm Std}}
\def\cstd{{\bf Std}}

\def\convol{{*}}

\def\park{{\bf a}}

\def\<{\langle}
\def\>{\rangle}

\def\C{\operatorname{\mathbb C}}
\def\K{\operatorname{\mathbb C}}
\def\Z{\operatorname{\mathbb Z}}
\def\N{\operatorname{\mathbb N}}

\def\gl{{\mathfrak gl}}
\def\glchap{\widehat\gl}

\def\F{{\bf F}}

\def\G{{\bf G}}

\def\V{{\bf V}}

\def\SG{{\mathfrak S}}

\def\dim{{\rm dim}}

\def\PF{{\rm PF}}
\def\PPF{{\rm PPF}}

\def\bij{{\phi}}
\newcommand{\free}[1]{\langle#1\rangle}
\def\conn{{\mathcal C}}
\def\L{{\mathcal L}}
\def\npn{{\bf n}}

\def\finer{\leq}
\def\MR{{\rm MR}}

\title[Free quasi-symmetric functions of level $l$]{Free quasi-symmetric functions 
of arbitrary level}

\author[J.-C.~Novelli and J.-Y.~Thibon]%
{Jean-Christophe Novelli and Jean-Yves Thibon}

\address[] {Institut Gaspard Monge, Universit\'e de Marne-la-Vall\'ee \\
5 Boulevard Descartes \\Champs-sur-Marne \\77454 Marne-la-Vall\'ee cedex 2 \\
FRANCE}
\email[Jean-Christophe Novelli]{novelli@univ-mlv.fr}
\email[Jean-Yves Thibon]{jyt@univ-mlv.fr} 
\date{}

\begin{document}

\begin{abstract}
We introduce analogues of the Hopf algebra of Free quasi-symmetric functions
with bases labelled by colored permutations.
As applications, we recover in a simple way the descent algebras associated
with wreath products $\Gamma\wr\SG_n$ and the corresponding generalizations of
quasi-symmetric functions. Finally, we obtain Hopf algebras of colored parking
functions, colored non-crossing partitions and parking functions of type $B$.
\end{abstract}

\maketitle

%%%%%%%%%%%%%%%%%%%%%%%%%%%%%%%%%%%%%%%%%%%%%%%%%%%%%%%%%%%%%%%%%%%%%%%%%%%%%%%
%%%%%%%%%%%%%%%%%%%%%%%%%%%%%%%%%%%%%%%%%%%%%%%%%%%%%%%%%%%%%%%%%%%%%%%%%%%%%%%
%%%%%%%%%%%%%%%%%%%%%%%%%%%%%%%%%%%%%%%%%%%%%%%%%%%%%%%%%%%%%%%%%%%%%%%%%%%%%%%
\section{Introduction}

The Hopf algebra of Free quasi-symmetric functions $\FQSym$ \cite{NCSF6}
is a certain algebra of noncommutative polynomials associated with the
sequence $(\SG_n)_{n\geq0}$ of all symmetric groups.
It is connected by Hopf homomorphisms to several other important algebras
associated with the same sequence of groups : Free symmetric functions (or
coplactic algebra) $\FSym$ \cite{PR,NCSF6}, Non-commutative symmetric
functions (or descent algebras) $\NCSF$ \cite{NCSF1}, Quasi-Symmetric
functions $\QSym$ \cite{Ge84}, Symmetric functions $\sym$, and also, Planar
binary trees $\PBT$ \cite{LR1,HNT}, Matrix quasi-symmetric functions $\MQSym$
\cite{NCSF6,Hiv}, Parking functions $\PQSym$ \cite{KW,NT}, and so on.

Among the many possible interpretations of $\sym$, we can mention the
identification as the representation ring of the tower of algebras
\begin{equation}
\C \to \C\SG_1 \to \C\SG_2 \to \cdots \to \C\SG_n \to \cdots,
\end{equation}
that is
\begin{equation}
\sym \simeq \oplus_{n\geq0} R(\C\SG_n),
\end{equation}
where $R(\C\SG_n)$ is the vector space spanned by isomorphism classes of
irreducible representations of $\SG_n$, the ring operations being induced by
direct sum and outer tensor product of representations \cite{Mcd}.

Another important interpretation of $\sym$ is as the support
of Fock space representations of
various infinite dimensional Lie algebras, in particular as the level $1$
irreducible highest weight representations of $\glchap_\infty$ (the infinite
rank Kac-Moody algebra of type $A_\infty$, with Dynkin diagram $\Z$,
see~\cite{Kac}).
 
The analogous level $l$ representations of this algebra can also be naturally
realized with products of $l$ copies of $\sym$, or as symmetric functions in
$l$ independent sets of variables
\begin{equation}
(\sym)^{\otimes l} \simeq \sym(X_0 ; \ldots ; X_{l-1}) =: \sym^{(l)},
\end{equation}
and these algebras are themselves the representation rings of wreath product
towers $(\Gamma\wr\SG_n)_{n\geq0}$, $\Gamma$ being a group with $l$ conjugacy
classes \cite{Mcd} (see also \cite{Zel,Wang}).

We shall therefore call for short $\sym(X_0 ; \ldots ; X_{l-1})$ the algebra
of symmetric functions of level $l$.
Our aim is to associate with $\sym^{(l)}$ analogues of the various Hopf
algebras mentionned at the beginning of this introduction.

We shall start with a level $l$ analogue of $\FQSym$, whose bases are labelled
by $l$-colored permutations. Imitating the embedding of $\NCSF$ in $\FQSym$,
we obtain a Hopf subalgebra of level $l$ called $\NCSF^{(l)}$, which turns out
to be dual to Poirier's quasi-symmetric functions, and whose homogenous
components can be endowed with an internal product, providing an analogue
of Solomon's descent algebras for wreath products.

The Mantaci-Reutenauer descent algebra arises as a natural Hopf subalgebra of
$\NCSF^{(l)}$ and its dual is computed in a straightforward way by means of an
appropriate Cauchy formula.

Finally, we introduce a Hopf algebra of colored parking functions,
and use it to define Hopf algebras structures on parking functions
and non-crossing partitions of type $B$.

%%%%%%%%%%%%%%%%%%%%%%%%%%%%%%%%%%%%%%%%%%%%%%%%%%%%%%%%%%%%%%%%%%%%%%%%%%%%%%%
{\it Acknowledgements} This research has been partially supported
by EC's IHRP Programme, grant HPRN-CT-2001-00272, ``Algebraic Combinatorics
in Europe".
%%%%%%%%%%%%%%%%%%%%%%%%%%%%%%%%%%%%%%%%%%%%%%%%%%%%%%%%%%%%%%%%%%%%%%%%%%%%%%%
%%%%%%%%%%%%%%%%%%%%%%%%%%%%%%%%%%%%%%%%%%%%%%%%%%%%%%%%%%%%%%%%%%%%%%%%%%%%%%%
\section{Free quasi-symmetric functions of level $l$}

\subsection{$l$-colored standardization}

We shall start with an $l$-colored alphabet
\begin{equation}
A = A^0 \sqcup A^1 \sqcup \cdots \sqcup A^{l-1},
\end{equation}
such that all $A^i$ are of the same cardinality $N$, which will be assumed to
be infinite in the sequel.
Let $C$ be the alphabet $\{c_0,\ldots,c_{l-1}\}$ and $B$ be the auxiliary
ordered alphabet $\{1,2,\ldots,N\}$ (the letter $C$ stands for \emph{colors}
and $B$ for \emph{basic}) so that $A$ can be identified to the cartesian
product $B\times C$:
\begin{equation}
A \simeq B \times C = \{ (b,c), b\in B,\ c\in C \}.
\end{equation}
Let $w$ be a word in $A$, represented as $(v,u)$ %$\bimot{u}{v}$
with $v\in B^*$ and $u\in C^*$. Then the \emph{colored standardized word}
$\cstd(w)$ of $w$ is
\begin{equation}
\cstd(w) := (\Std(v),u), %\bimot{u}{\Std(v)}
\end{equation}
where $\Std(v)$ is the usual standardization on words.

Recall that the standardization process sends a word $w$ of length $n$ to a
permutation $\Std(w)\in\SG_n$ called the \emph{standardized} of $w$ defined as
the permutation obtained by iteratively scanning $w$ from left to right, and
labelling $1,2,\ldots$ the occurrences of its smallest letter, then numbering
the occurrences of the next one, and so on. Alternatively, $\Std(w)$ is the
permutation having the same inversions as $w$.

%For example, $\Std({abcadbcaa})=157296834$:
%\begin{equation}
%\begin{array}{ccccccccc}
%a  & b  & c  & a  & d  & b  & c  & a  & a\\
%a_1& b_5& c_7& a_2& d_9& b_6& c_8& a_3& a_4\\
%1   &5   &7   &2   &9   &6   &8   &3   &4
%\end{array}
%\end{equation}

%%%%%%%%%%%%%%%%%%%%%%%%%%%%%%%%%%%%%%%%%%%%%%%%%%%
\subsection{$\FQSym^{(l)}$ and $\FQSym^{(\Gamma)}$}

A \emph{colored permutation} is a pair $(\sigma,u)$, with $\sigma\in\SG_n$ and
$u\in C^n$, the integer $n$ being the \emph{size} of this permutation.

\begin{definition}
The \emph{dual free $l$-quasi-ribbon} $\G_{\sigma,u}$ labelled by a
colored permutation $(\sigma,u)$ of size $n$ is the noncommutative
polynomial
\begin{equation}
\G_{\sigma,u} := \sum_{w\in A^n ; \cstd(w)=(\sigma,u)} w \quad\in\Z\free{A}.
\end{equation}
\end{definition}

Recall that the \emph{convolution} of two permutations $\sigma$ and $\mu$ is
the set $\sigma\convol\mu$ (identified with the formal sum of its elements)
of permutations $\tau$ such that the standardized word of the $|\sigma|$ first
letters of $\tau$ is $\sigma$ and the standardized word of the remaining
letters of $\tau$ is $\mu$ (see~\cite{Reu}).

%For an element $v=(v_1,v_2,\ldots)$ of $B$ and an integer $k$, denote by
%$v[k]$ the \emph{shifted word} $(v_1+k)\cdots (v_n+k)$, e.g., $312[4]=756$.
%The \emph{shifted concatenation} of two words $v$ and $v'$ is defined by
%%
%\begin{equation}
%v\bullet v'= v\cdot v'[k]
%\end{equation}
%%
%where $k$ is the length of $v$.

\begin{theorem}
\label{prodG}
Let $(\sigma',u')$ and $(\sigma'',u'')$ be colored permutations.
Then
\begin{equation}
\G_{\sigma',u'}\,\,\G_{\sigma'',u''} = \sum_{\sigma\in \sigma'\convol\sigma''}
\G_{\sigma,u'\cdot u''},
\end{equation}
where $w_1\cdot w_2$ is the word obtained by concatenating $w_1$ and $w_2$.
Therefore, the dual free $l$-quasi-ribbons span a $\Z$-subalgebra of the free
associative algebra.
\end{theorem}

Moreover, one defines a coproduct on the $\G$ functions by

\begin{equation}
\label{deltaG}
\Delta \G_{\sigma,u} := \sum_{i=0}^n
\G_{(\sigma,u)_{[1,i]}} \otimes \G_{(\sigma,u)_{[i+1,n]}},
\end{equation}
where $n$ is the size of $\sigma$ and $(\sigma,u)_{[a,b]}$ is the standardized
colored permutation of the pair $(\sigma',u')$ where $\sigma'$ is the subword
of $\sigma$ containing the letters of the interval $[a,b]$, and $u'$ the
corresponding subword of $u$.

For example,
\begin{equation}
\begin{split}
\Delta \G_{3142,2412} =
&\ 1\otimes \G_{3142,2412} + \G_{1,4}\otimes \G_{231,212} +
   \G_{12,42}\otimes \G_{12,21} \\
& + \G_{312,242}\otimes \G_{1,1} + \G_{3142,2412}\otimes 1.
\end{split}
\end{equation}

\begin{theorem}
The coproduct is an algebra homomorphism, so that $\FQSym^{(l)}$ is a
graded bialgebra. Moreover, it is a Hopf algebra.
\end{theorem}

\begin{definition}
The \emph{free $l$-quasi-ribbon} $\F_{\sigma,u}$ labelled by a colored
permutation $(\sigma,u)$ is the noncommutative polynomial
\begin{equation}
\F_{\sigma,u} := \G_{\sigma^{-1},u\cdot\sigma^{-1}},
\end{equation}
where the action of a permutation on the right of a word permutes the
positions of the letters of the word.
\end{definition}

For example,
\begin{equation}
\F_{3142,2142} = \G_{2413,1422}\,.
\end{equation}

The product and coproduct of the $\F_{\sigma,u}$ can be easily described in
terms of shifted shuffle and deconcatenation of colored permutations.

Let us define a scalar product on $\FQSym^{(l)}$ by
\begin{equation}
\langle \F_{\sigma,u} , \G_{\sigma',u'} \rangle :=
  \delta_{ \sigma,\sigma'} \delta_{u,u'},
\end{equation}
where $\delta$ is the Kronecker symbol.

\begin{theorem}
For any $U,V,W\in\FQSym^{(l)}$,
\begin{equation}
\langle \Delta U, V\otimes W \rangle =
\langle U, V W \rangle,
\end{equation}
so that, $\FQSym^{(l)}$ is self-dual: the map
$\F_{\sigma,u} \mapsto {\G_{\sigma,u}}^*$ is an isomorphism from
$\FQSym^{(l)}$ to its graded dual.
\end{theorem}

\begin{note}{\rm
Let $\bij$ be any bijection from $C$ to $C$, extended to words by
concatenation. Then if one defines the free $l$-quasi-ribbon as
\begin{equation}
\F_{\sigma,u} := \G_{\sigma^{-1},\bij(u)\cdot\sigma^{-1}},
\end{equation}
the previous theorems remain valid since one only permutes the labels of
the basis $(\F_{\sigma,u})$.

Moreover, if $C$ has a group structure, the colored permutations
$(\sigma,u)\in\SG_n\times C^n$ can be interpreted as elements of the
semi-direct product $H_n := \SG_n\ltimes C^n$ with multiplication rule
\begin{equation}
(\sigma ; c_1,\ldots,c_n) \cdot (\tau ; d_1,\ldots,d_n) :=
(\sigma\tau ; c_{\tau(1)}d_1, \ldots, c_{\tau(n)}d_n).
\end{equation}
In this case, one can choose $\bij(\gamma):=\gamma^{-1}$ and define the scalar
product as before, so that the adjoint basis of the $(\G_{h})$ becomes
$\F_h := \G_{h^{-1}}$.
In the sequel, we will be mainly interested in the case $C:=\Z/l\Z$, and we
will indeed make that choice for $\bij$.
}
\end{note}

%$V=\oplus_{\gamma\in\Gamma} V_\gamma$
%matrices $\gamma$-circulantes par blocs -> $\FQSym^{(\Gamma)}$.

%%%%%%%%%%%%%%%%%%%%%%%%%%%%%%%%
\subsection{Algebraic structure}

Recall that a permutation $\sigma$ of size $n$ is \emph{connected}
\cite{MR,NCSF6} if, for any $i<n$, the set $\{\sigma(1),\ldots,\sigma(i)\}$ is
different from $\{1,\ldots,i\}$.

We denote by $\conn$ the set of connected permutations, and by
$c_n:=|\conn_n|$ the number of such permutations in $\SG_n$. For later
reference, we recall that the generating series of $c_n$ is
\begin{equation*}
c(t) := \sum_{n\ge 1} c_n t^n
  = 1 - \left(\sum_{n\ge 0} n! t^n\right)^{-1}\\
  = t+{t}^{2}+3\,{t}^{3}+13\,{t}^{4}+71\,{t}^{5}+461\,{t}^{6} +O(t^7)\,.
\end{equation*}

Let the \emph{connected colored permutations} be the $(\sigma,u)$ with
$\sigma$ connected and $u$ arbitrary. Their generating series is given by
$c(lt)$.

It follows from \cite{NCSF6} that $\FQSym^{(l)}$ is free over the set
$\G_{\sigma,u}$ (or $\F_{\sigma,u}$), where $(\sigma,u)$ is connected.

Since $\FQSym^{(l)}$ is self-dual, it is also cofree.

%%%%%%%%%%%%%%%%%%%%%%%%%%%%%%%
\subsection{Primitive elements}

Let $\L^{(l)}$ be the primitive Lie algebra of $\FQSym^{(l)}$.
Since $\Delta$ is not cocommutative, $\FQSym^{(l)}$ cannot be the universal
enveloping algebra of $\L^{(l)}$.
But since it is cofree, it is, according to~\cite{LRdip}, the universal
enveloping dipterous algebra of its primitive part $\L^{(l)}$.
Let $d_n = \dim\, \L^{(l)}_n$.

Recall that the \emph{shifted concatenation} $w\bullet w'$ of two elements $w$
and $w'$ of $\N^*$, is the word obtained by concatenating to $w$ the word
obtained by shifting all letters of $w'$ by the length of $w$. We extend
it to colored permutations by simply concatenating the colors and
concatenating \emph{with shift} the permutations.
Let $\G^{\sigma,u}$ be the multiplicative basis defined by
$\G^{\sigma,u}=\G_{\sigma_1,u_1}\cdots\G_{\sigma_r,u_r}$ where
$(\sigma,u)=(\sigma_1,u_1)\bullet\cdots\bullet(\sigma_r,u_r)$ is the unique
maximal factorization of $(\sigma,u)\in\SG_n\times C^n$ into connected colored
permutations.

\begin{proposition}
Let $\V_{\sigma,u}$ be the adjoint basis of $\G^{\sigma,u}$.
Then, the family $(\V_{\alpha,u})_{\alpha\in\conn}$ is a basis of $\L^{(l)}$.
In particular, we have $d_n=l^n c_n$.
\end{proposition}

As in~\cite{NCSF6}, we conjecture that $\L^{(l)}$ is free.

%%%%%%%%%%%%%%%%%%%%%%%%%%%%%%%%%%%%%%%%%%%%%%%%%%%%%%%%%%%%%%%%%%%%%%%%%%%%%%%
%%%%%%%%%%%%%%%%%%%%%%%%%%%%%%%%%%%%%%%%%%%%%%%%%%%%%%%%%%%%%%%%%%%%%%%%%%%%%%%
%%%%%%%%%%%%%%%%%%%%%%%%%%%%%%%%%%%%%%%%%%%%%%%%%%%%%%%%%%%%%%%%%%%%%%%%%%%%%%%
%%%%%%%%%%%%%%%%%%%%%%%%%%%%%%%%%%%%%%%%%%%%%%%%%%%%%%%%%%%%%%%%%%%%%%%%%%%%%%%
\section{Non-commutative symmetric functions of level $l$}

Following McMahon~\cite{McM}, we define an \emph{$l$-partite number} $\npn$
as a column vector in $\N^l$, and a \emph{vector composition of $\npn$} of
weight $|\npn|:=\sum_{1}^l{n_i}$ and length $m$ as a $l\times m$ matrix
$\bf I$ of nonnegative integers, with row sums vector $\npn$ and no zero
column.

For example,
\begin{equation}
\label{exM}
{\bf I} =
\begin{pmatrix}
1 & 0 & 2 & 1 \\
0 & 3 & 1 & 1 \\
4 & 2 & 1 & 3 \\
\end{pmatrix}
\end{equation}
is a vector composition (or a $3$-composition, for short)
of the $3$-partite number
\scriptsize$\begin{pmatrix} 4\\ 5\\ 10\end{pmatrix}$ \normalsize
of weight $19$ and length $4$.

For each $\npn\in\N^l$  of weight $|\npn|=n$, we define a level $l$
\emph{complete homogeneous symmetric function} as
\begin{equation}
S_{\npn} := \sum_{u ; |u|_i=n_i} \G_{1\cdots n, u}.
\end{equation}
It is the sum of all possible colorings of the identity permutation with $n_i$
occurrences of color $i$ for each $i$.

Let $\NCSF^{(l)}$ be the subalgebra of $\FQSym^{(l)}$ generated by the
$S_{\npn}$ (with the convention $S_{\bf 0}=1$).
The Hilbert series of $\NCSF^{(l)}$ is easily found to be
\begin{equation}
S_l(t) := \sum_{n} {\dim\ \NCSF_n^{(l)}t^n} = \frac{(1-t)^l}{2(1-t)^l-1}.
\end{equation}

\begin{theorem}
$\NCSF^{(l)}$  is free over the set $\{S_{\npn}, |\npn|>0 \}$.
Moreover, $\NCSF^{(l)}$ is a Hopf subalgebra of $\FQSym^{(l)}$.

The coproduct of the generators is given by
\begin{equation}
\Delta S_\npn = \sum_{{\bf i}+{\bf j}= {\bf n}} S_{\bf i}\otimes S_{\bf j},
\end{equation}
where the sum ${\bf i}+{\bf j}$ is taken in the space $\N^l$. In particular,
$\NCSF^{(l)}$ is cocommutative.
\end{theorem}

%For example,

We can therefore introduce the basis of products of level $l$ complete
function, labelled by $l$-compositions
\begin{equation}
S^{\bf I} = S_{{\bf i}_1} \cdots S_{{\bf i}_m},
\end{equation}
where ${\bf i}_1,\cdots,{\bf i}_m$ are the columns of $\bf I$.

\begin{theorem}
If $C$ has a group structure, $\NCSF_n^{(l)}$ becomes a subalgebra of
$\C[C\wr\SG_n]$ under the identification $\G_h \mapsto h$.
\end{theorem}

This provides an analogue of Solomon's descent algebra for the wreath
product $C\wr\SG_n$. The proof amounts to check that the Patras descent
algebra of a graded bialgebra \cite{Patras} can be adapted to $\N^l$-graded
bialgebras.

As in the case $l=1$, we define the \emph{internal product} $*$ as being
opposite to the law induced by the group algebra. It can be computed by the
following splitting formula, which is a straightforward generalization of the
level 1 version.

\begin{proposition}
Let $\mu_r: (\NCSF^{(l)})^{\otimes r} \to \NCSF^{(l)}$ be the product map.
Let $\Delta^{(r)} : (\NCSF^{(l)}) \to (\NCSF^{(l)})^{\otimes r}$ be the
$r$-fold coproduct, and $*_r$ be the extension of the internal product to
$(\NCSF^{(l)})^{\otimes r}$.
Then, for $F_1,\ldots,F_r$, and $G\in\NCSF^{(l)}$,

\begin{equation}
(F_1\cdots F_r) * G = \mu_r [ (F_1\otimes\cdots\otimes F_r) *_r
\Delta^{(r)}G ].
\end{equation}
\end{proposition}

The group law of $C$ is needed only for the evaluation of the
product of one-part complete functions $S_{\bf m}*S_{\bf n}$.

\begin{example}
With $l=2$ and $C=\Z/2\Z$,

\scriptsize
\begin{equation*}
\begin{split}
S^{\hbox{\scriptsize$\begin{pmatrix}1&0\\1&1\end{pmatrix}$}} *
S^{\hbox{\scriptsize$\begin{pmatrix}0&2\\1&0\end{pmatrix}$}} &=
\mu_2 \left[ \left(
    S^{\hbox{\scriptsize$\begin{pmatrix}1\\1\end{pmatrix}$}}\otimes
      S^{\hbox{\scriptsize$\begin{pmatrix}0\\1\end{pmatrix}$}} \right)
    *_2
    \Delta S^{\hbox{\scriptsize$\begin{pmatrix}0&2\\1&0\end{pmatrix}$}}
    \right]\\
& =  \left(S^{\hbox{\scriptsize$\begin{pmatrix}1\\1\end{pmatrix}$}}*
           S^{\hbox{\scriptsize$\begin{pmatrix}2\\0\end{pmatrix}$}}\right)
     \left(S^{\hbox{\scriptsize$\begin{pmatrix}0\\1\end{pmatrix}$}}*
           S^{\hbox{\scriptsize$\begin{pmatrix}1\\0\end{pmatrix}$}}\right)
  +  \left(S^{\hbox{\scriptsize$\begin{pmatrix}1\\1\end{pmatrix}$}}*
           S^{\hbox{\scriptsize$\begin{pmatrix}0&1\\1&0\end{pmatrix}$}}\right)
     \left(S^{\hbox{\scriptsize$\begin{pmatrix}0\\1\end{pmatrix}$}}*
           S^{\hbox{\scriptsize$\begin{pmatrix}1\\0\end{pmatrix}$}}\right) \\
& = S^{\hbox{\scriptsize$\begin{pmatrix}1&1\\1&0\end{pmatrix}$}} 
  + S^{\hbox{\scriptsize$\begin{pmatrix}1&1&0\\0&0&1\end{pmatrix}$}} 
  + S^{\hbox{\scriptsize$\begin{pmatrix}0&0&0\\1&1&1\end{pmatrix}$}}.
\end{split}
\end{equation*}
\end{example}

Recall that the underlying colored alphabet $A$ can be seen as
$A^0 \sqcup \cdots \sqcup A^{l-1}$, with $A^i = \{ a^{(i)}_j, j\geq1 \}$.
Let ${\bf x} = (x^{(0)}, \ldots, x^{(l-1)})$, where the $x^{(i)}$ are $l$
commuting variables.
In terms of $A$, the generating function of the complete functions can be
written as
\begin{equation}
\sigma_{\bf x}(A) = \prod_{i\geq0}^{\rightarrow}
\left(1-\sum_{0\leq j\leq l-1} x^{(j)} a_{i}^{(j)} \right)^{-1}
= \sum_\npn {S_{\bf n}(A) {\bf x}^{\bf n}},
\end{equation}
where ${\bf x}^{\bf n} = (x^{(0)})^{n_0} \cdots (x^{(l-1)})^{n_{l-1}}$.

This realization gives rise to a Cauchy formula (see~\cite{NCSF2} for the
$l=1$ case), which in turn allows one to identify the dual of $\NCSF^{(l)}$ with
an algebra introduced by S. Poirier in~\cite{Poi}.

%%%%%%%%%%%%%%%%%%%%%%%%%%%%%%%%%%%%%%%%%%%%%%%%%%%%%%%%%%%%%%%%%%%%%%%%%%%%%%%
%%%%%%%%%%%%%%%%%%%%%%%%%%%%%%%%%%%%%%%%%%%%%%%%%%%%%%%%%%%%%%%%%%%%%%%%%%%%%%%
\section{Quasi-symmetric functions of level $l$}

%In this section, we shall assume that $C$ is the cyclic group $\Z/l\Z$ so that
%$\F_{\sigma,u}$ is defined as $\G_{\sigma^{-1},u'\cdot\sigma^{-1}}$,
%where $u'$ is the opposite of $u$ in the abelian group $(\Z/l\Z)^n$.

%%%%%%%%%%%%%%%%%%%%%%%%%%%%%%%%%%%%%%%%
\subsection{Cauchy formula of level $l$}

Let $X = X^0 \sqcup \cdots \sqcup X^{l-1}$, where $X^i=\{ x_j^{(i)},j\geq1\}$
be an $l$-colored alphabet of commutative variables, also commuting with $A$.
Imitating the level $1$ case (see~\cite{NCSF6}), we define the Cauchy kernel

\begin{equation}
K(X,A) = \prod_{j\geq1}^{\rightarrow}
\sigma_{\left(x_j^{(0)}, \ldots, x_j^{(l-1)}\right)} (A).
\end{equation}

Expanding on the basis $S^{\bf I}$ of $\NCSF^{(l)}$, we get as coefficients
what can be called the \emph{level $l$ monomial quasi-symmetric functions}
$M_{\bf I}(X)$

\begin{equation}
K(X,A) = \sum_{\bf I} M_{\bf I}(X) S^{\bf I}(A),
\end{equation}
defined by
\begin{equation}
M_{\bf I}(X) = \sum_{j_1<\cdots<j_m}
{\bf x}^{{\bf i}_1}_{j_1} \cdots
{\bf x}^{{\bf i}_m}_{j_m},
\end{equation}
with ${\bf I}=({\bf i}_1,\ldots,{\bf i}_m)$.

These last functions form a basis of a subalgebra $\QSym^{(l)}$ of $\K[X]$,
which we shall call the \emph{algebra of quasi-symmetric functions of level
$l$}.

%%%%%%%%%%%%%%%%%%%%%%%%%%%%%%%%%%%%%%%%%%%%%%%%
\subsection{Poirier's Quasi-symmetric functions}

The functions $M_{\bf I}(X)$ can be recognized as a basis of one of the
algebras introduced by Poirier: the $M_{\bf I}$ coincide with the $M_{(C,v)}$
defined in~\cite{Poi}, p.~324, formula (1), up to indexation.

Following Poirier, we introduce the level $l$ quasi-ribbon functions by
summing over an order on $l$-compositions:
an $l$-composition $C$ is finer than $C'$, and we write $C\finer C'$, if
$C'$ can be obtained by repeatedly summing up two consecutive columns of $C$
such that no non-zero element of the left one is strictly below a non-zero
element of the right one.
%any two columns summed together (a left one and a right one) have no non-zero
%element

This order can be described in a much easier and natural way if one recodes
an $l$-composition ${\bf I}$ as a pair of words, the first one $d({\bf I})$
being the set of sums of the elements of the first $k$ columns of $\bf I$, the
second one $c({\bf I})$ being obtained by concatenating the words
$i^{{\bf I}_{i,j}}$ while reading of $\bf I$ by columns, from top to bottom
and from left to right.
For example, the $3$-composition of Equation~(\ref{exM}) satisfies
\begin{equation}
d({\bf I}) = \{5, 10, 14, 19\} \text{\quad and\quad}
c({\bf I}) = 13333\, 22233\, 1123\, 12333\,.
\end{equation}
Moreover, this recoding is a bijection if the two words $d({\bf I})$ and
$c({\bf I})$ are such that the descent set of $c({\bf I})$ is a subset of
$d({\bf I})$.
The order previously defined on $l$-compositions is in this context the
inclusion order on sets $d$: $(d',c)\finer (d,c)$ iff $d'\subseteq d$.

It allows us to define the \emph{level $l$ quasi-ribbon functions} $F_{\bf I}$
by
\begin{equation}
F_{\bf I} = \sum_{{\bf I'}\finer {\bf I}} M_{\bf I'}.
\end{equation}
Notice that this last description of the order $\finer$ is reminiscent of the
order $\finer'$ on descent sets used in the context of quasi-symmetric
functions and non-commutative symmetric functions: more precisely, since it
does not modify the word $c({\bf I})$, the order $\finer$ restricted to
$l$-compositions of weight $n$ amounts for $l^n$ copies of the order
$\finer'$.
The computation of its M\"obius function is therefore straightforward.

Moreover, one can directly obtain the $F_{\bf I}$ as the commutative image
of certain $\F_{\sigma,u}$: any pair $(\sigma,u)$ such that $\sigma$ has
descent set $d({\bf I})$ and $u=c({\bf I})$ will do.

%%%%%%%%%%%%%%%%%%%%%%%%%%%%%%%%%%%%%%%%%%%%%%%%%%%%%%%%%%%%%%%%%%%%%%%%%%%%%%%
%%%%%%%%%%%%%%%%%%%%%%%%%%%%%%%%%%%%%%%%%%%%%%%%%%%%%%%%%%%%%%%%%%%%%%%%%%%%%%%
%%%%%%%%%%%%%%%%%%%%%%%%%%%%%%%%%%%%%%%%%%%%%%%%%%%%%%%%%%%%%%%%%%%%%%%%%%%%%%%
\section{The Mantaci-Reutenauer algebra}

Let ${\bf e}_i$ be the canonical basis of $\N^l$. For $n\geq1$, let
\begin{equation}
S_n^{(i)} = S_{n\cdot {\bf e}_i} \in\NCSF^{(l)},
\end{equation}
be the \emph{monochromatic complete symmetric functions}.

\begin{proposition}
The $S_n^{(i)}$ generate a Hopf-subalgebra $\MR^{(l)}$ of $\NCSF^{(l)}$, which
is isomorphic to the Mantaci-Reutenauer descent algebra of the wreath products
$\SG_n^{(l)} = (\Z/l\Z) \wr\SG_n$.
\end{proposition}

It follows in particular that $\MR^{(l)}$ is stable under the composition
product induced by the group structure of $\SG_n^{(l)}$.
The bases of $\MR^{(l)}$ are labelled by colored compositions (see below).

The duality is easily worked out by means of the appropriate Cauchy kernel.
The generating function of the complete functions is
\begin{equation}
\sigma^{\MR}_{\bf x}(A) := 1 + \sum_{j=0}^{l-1} \sum_{n\geq1}
S_n^{(j)}.(x^{(j)})^n,
\end{equation}
and the Cauchy kernel is as usual
\begin{equation}
K^{\MR}(X,A) = \prod_{i\geq1}^\rightarrow \sigma^{\MR}_{{\bf x}_i}(A)
= \sum_{(I,u)} M_{(I,u)}(X) S^{(I,u)}(A),
\end{equation}
where $(I,u)$ runs over colored compositions
$(I,u) = ((i_1,\ldots,i_m),(u_1,\ldots,u_m))$ that is, pairs formed with a
composition and a color vector of the same length. The $M_{I,u}$ are
called the \emph{monochromatic monomial quasi-symmetric functions} and satisfy
\begin{equation}
M_{(I,u)}(X) = \sum_{j_1<\cdots<j_m}
(x_{j_1}^{(u_1)})^{i_1} \cdots (x_{j_m}^{(u_m)})^{i_m}.
\end{equation}

\begin{proposition}
The $M_{(I,u)}$ span a subalgebra of $\C[X]$ which can be identified with the
graded dual of $\MR^{(l)}$ through the pairing
\begin{equation}
\langle M_{(I,u)}, S^{(J,v)} \rangle = \delta_{I,J} \delta_{u,v},
\end{equation}
where $\delta$ is the Kronecker symbol.
\end{proposition}

Note that this algebra can also be obtained by assuming the relations
\begin{equation}
x_i^{(p)} x_i^{(q)} = 0, \text{\ for $p\not=q$}
\end{equation}
on the variables of $\QSym^{(l)}$.

Baumann and Hohlweg have another construction of the dual of $\MR^{(l)}$
\cite{BH} (implicitly defined in~\cite{Poi}, Lemma~11).

%%%%%%%%%%%%%%%%%%%%%%%%%%%%%%%%%%%%%%%%%%%%%%%%%%%%%%%%%%%%%%%%%%%%%%%%%%%%%%%
%%%%%%%%%%%%%%%%%%%%%%%%%%%%%%%%%%%%%%%%%%%%%%%%%%%%%%%%%%%%%%%%%%%%%%%%%%%%%%%
%%%%%%%%%%%%%%%%%%%%%%%%%%%%%%%%%%%%%%%%%%%%%%%%%%%%%%%%%%%%%%%%%%%%%%%%%%%%%%%
\section{Level $l$ parking quasi-symmetric functions}

\subsection{Usual parking functions}

Recall that a \emph{parking function} on $[n]=\{1,2,\ldots,n\}$ is a word
$\park=a_1a_2\cdots a_n$ of length $n$ on $[n]$ whose nondecreasing
rearrangement $\park^\uparrow=a'_1a'_2\cdots a'_n$ satisfies $a'_i\le i$ for
all $i$.
Let $\PF_n$ be the set of such words.
It is well-known that $|\PF_n|=(n+1)^{n-1}$.

Gessel introduced in 1997 (see~\cite{Stan2}) the notion of \emph{prime parking
function}. One says that $\park$ has a \emph{breakpoint} at $b$ if
$|\{\park_i\le b\}|=b$. The set of all breakpoints of $\park$ is denoted by
$\bp(\park)$.
Then, $\park\in \PF_n$ is prime if $\bp(\park)=\{n\}$.

Let $\PPF_n\subset\PF_n$ be the set of prime parking functions on $[n]$.
It can easily be shown that $|\PPF_n|=(n-1)^{n-1}$ (see~\cite{Stan2}).

We will finally need one last notion: $\park$ has a \emph{match} at $b$ if
$|\{\park_i< b\}|=b-1$ and $|\{\park_i\leq b\}|\geq b$. The set of all matches
of $\park$ is denoted by $\match(\park)$.

We will now define generalizations of the usual parking functions to any level
in such a way that they build up a Hopf algebra in the same way as
in~\cite{NT}.

\subsection{Colored parking functions}

Let $l$ be an integer, representing the number of allowed colors.
A \emph{colored parking function} of level $l$ and size $n$ is a pair
composed of a parking function of length $n$ and a word on $[l]$ of length
$l$.

Since there is no restriction on the coloring, it is obvious that there are
$l^n (n+1)^{n-1}$ colored parking functions of level $l$ and size $n$.

It is known that the convolution of two parking functions contains only
parking functions, so one easily builds as in~\cite{NT} an algebra
$\PQSym^{(l)}$ indexed by colored parking functions:

\begin{equation}
\G_{(\park',u')} \G_{(\park'',u'')} = \sum_{\park\in \park'\convol\park''}
  \G_{(\park,u'\cdot  u'')}.
\end{equation}

Moreover, one defines a coproduct on the $\G$ functions by
\begin{equation}
\Delta\G_{(\park,u)} = \sum_{i\in \bp(\park)} 
\G_{(\park,u)_{[1,i]}} \otimes \G_{(\park,u)_{[i+1,n]}}
\end{equation}
where $n$ is the size of $\park$ and $(\park,u)_{[a,b]}$ is the parkized
colored parking function of the pair $(\park',u')$ where $\park'$ is the
subword of $\park$ containing the letters of the interval $[a,b]$, and $u'$
the corresponding subword of $u$.

\begin{theorem}
The coproduct is an algebra homomorphism, so that $\PQSym^{(l)}$ is a
graded bialgebra. Moreover, it is a Hopf algebra.
\end{theorem}

\subsection{Parking functions of type $B$}

In~\cite{Rei}, Reiner defined non-crossing partitions of type $B$ by
analogy to the classical case. In our context, he defined the level $2$ case.
It allowed him to derive, by analogy with a simple representation theoretical
result, a definition of parking functions of type $B$ as the words on $[n]$ of
size $n$.

We shall build another set of
words, also enumerated by $n^n$ that sheds light on the relation between
type-$A$ and type-$B$ parking functions and provides a natural Hopf algebra
structure on the latter.

First, fix two colors $0<1$. We say that a pair of words $(\park,u)$ composed
of a parking function and a binary colored word is a
\emph{level $2$ parking function} if
\begin{itemize}
\item the only elements of $\park$ that can have color $1$ are the matches of
$\park$.
\item for all element of $\park$ of color $1$, all letters equal to it and to
its left also have color $1$,
\item all elements of $\park$ have at least once the color $0$.
\end{itemize}

For example, there are $27$ level $2$ parking functions of size $3$: there are
the $16$ usual ones all with full color $0$, and the eleven new elements
\begin{equation}
\begin{split}
& (111,100), (111,110), (112,100), (121,100), (211,010),\\
& (113,100), (131,100), (311,010),
(122,010), (212,100), (221,100). \\
\end{split}
\end{equation}

The first time the first rule applies is with $n=4$, where one has to discard
the words $(1122,0010)$ and $(1122,1010)$ since $2$ is not a match of $1122$.
On the other hand, both words $(1133,0010)$ and $(1133,1010)$ are
$B_4$-parking functions since $1$ and $3$ are matches of $1133$.

\begin{theorem}
The restriction of $\PQSym^{(2)}$ to the $\G$ elements indexed by level $2$
parking functions is a Hopf subalgebra of $\PQSym^{(2)}$.
\end{theorem}

\subsection{Non-crossing partitions of type $B$}

%In~\cite{Rei}, Reiner defined non-crossing partitions of type $B$ by
%analogy to the classical case (see also~\cite{Bi1,Stan}).
%In our context, he defined the level $2$ case.

%Let us extend his definition to any level. First, let us remark that in the
Remark that in the level $1$ case, non-crossing partitions are in bijection
with non-decreasing parking functions.
To extend this correspondence to type $B$, let us start with a non-decreasing
parking function ${\bf b}$ (with no color). We factor it into the maximal
shifted concatenation of prime non-decreasing parking functions, and we choose
a color, here 0 or 1, for each factor. We obtain in this way $\binom{2n}{n}$
words $\pi$, which can be identified with
\emph{type $B$ non-crossing partitions}.

Let
\begin{equation}
{\bf P}^{\pi}=\sum_{{\bf a}\sim{\pi}}\F_{\bf a}\,
\end{equation}
where $\sim$ denotes equality up to rearrangement of the letters. 
Then,
\begin{theorem}
The ${\bf P}^{\pi}$, where $\pi$ runs over the above set of non-decreasing
signed parking functions, form the basis of a cocommutative Hopf subalgebra
of $\PQSym^{(2)}$.
\end{theorem}

All this can be extended to higher levels in a straightforward way: allow each
prime non-decreasing parking function to choose any color among $l$ and use
the factorization as above. Since non-decreasing parking functions are in
bijection with Dyck words, the choice can be described as: each block of a
Dyck word with no return-to-zero, chooses one color among $l$. In this
version, the generating series is obviously given by
\begin{equation}
\frac{1}{1- l\frac{1-\sqrt{1-4t}}{2}}.
\end{equation}
For $l=3$, we obtain the sequence A007854 of~\cite{Slo}.

%%%%%%%%%%%%%%%%%%%%%%%%%%%%%%%%%%%%%%%%%%%%%%%%%%%%%%%%%%%%%%%%%%%%%%%%%%%%%%%
%%%%%%%%%%%%%%%%%%%%%%%%%%%%%%%%%%%%%%%%%%%%%%%%%%%%%%%%%%%%%%%%%%%%%%%%%%%%%%%
%%%%%%%%%%%%%%%%%%%%%%%%%%%%%%%%%%%%%%%%%%%%%%%%%%%%%%%%%%%%%%%%%%%%%%%%%%%%%%%

\end{document}